\documentclass[12pt]{amsart}
\usepackage[english]{babel}
\usepackage{graphicx}
\usepackage{framed}
\usepackage[normalem]{ulem}
\usepackage{amsmath}
\usepackage{amsthm}
\usepackage{amssymb}
\usepackage{amsfonts}
\usepackage{enumerate}
\usepackage{verbatim}
\usepackage[utf8]{inputenc}
\usepackage[top=1 in,bottom=1in, left=1 in, right=1 in]{geometry}

\newcommand{\R}{\mathbb{R}}

\newcommand{\C}{\mathbb{C}}
\newcommand{\Gm}{\mathbb{G}_m}

\newcommand{\Z}{\mathbb{Z}}
\newcommand{\Res}{{\rm Res}}

\renewcommand{\emptyset}{\varnothing}

\theoremstyle{theorem}
\newtheorem{theorem}{Theorem}
\newtheorem{conjecture}[theorem]{Conjecture}
\newtheorem{lemma}[theorem]{Lemma}

\newtheorem{prop}[theorem]{Proposition}
\newtheorem{corollary}[theorem]{Corollary}
\theoremstyle{definition}
\newtheorem{definition}[theorem]{Definition}
\newtheorem{remark}[theorem]{Remark}
\newtheorem{example}[theorem]{Example}

\newcommand{\lct}{\text{\rm lct}}
\newcommand{\Aut}{\text{\rm Aut}}

\newcommand{\SL}{\operatorname{SL}}
\newcommand{\GL}{\operatorname{GL}}
\newcommand{\rk}{\operatorname{rk}}
\newcommand{\symmrk}{\operatorname{symmrk}}
\newcommand{\val}{\operatorname{val}}
\newcommand{\id}{\text{id}}

\newcommand{\ord}{\operatorname{ord}}
\newcommand{\symm}{{D}}

\title[{G-stable Rank Of Symmetric Tensors And Log Canonical Threshold}]{G-stable Rank Of Symmetric Tensors\\ And Log Canonical Threshold}
\author{Zhi Jiang}
\date{}

\begin{document}

\maketitle
\begin{abstract}
    Shitov recently gave a counterexample to Comon's conjecture that the symmetric tensor rank and tensor rank of a symmetric tensor are the same. In this paper we show that an analog of Comon's conjecture for the $G$-stable rank introduced by Derksen is true: the symmetric $G$-stable rank and $G$-stable rank of a symmetric tensor are the same. We also show that the log-canonical threshold of a singularity is bounded by the $G$-stable rank of the defining ideal.
\end{abstract}

\section{Introduction}

An order $d$ tensor is a vector in a tensor product of $d$ vector spaces. 
The are several generalizations of the rank of a matrix  to tensors of order $\geq 3$, for example
the tensor rank, border rank, sub-rank, slice rank and $G$-stable rank. A simple tensor in $V_1\otimes V_2\otimes\cdots \otimes V_d$ is a tensor of the form
$v_1\otimes v_2\otimes \cdots\otimes v_d$, where $v_i\in V_i$. The
tensor rank of $T\in V_1\otimes V_2\otimes \cdots \otimes V_d$ is the smallest number of simple tensors
that sum up to $T$.
The $G$-stable rank of a tensor
was introduced by Derksen in \cite{G-stable rank}. The slice rank and $G$-stable rank have been used to find bounds for the cap set problem (see~\cite{G-stable rank}, \cite{EllenbergGijswijt}, \cite{Jiang}, \cite{Tao}).

If $V_1=V_2=\cdots=V_d=V$ then there is a natural
action of the symmetric group $S_d$ on $V^{\otimes d}=V\otimes V\otimes \cdots\otimes V$. A tensor
invariant under this action is called a symmetric tensor of order $d$.
The Waring rank or symmetric rank of a symmetric
tensor $T\in V^{\otimes d}$ is the smallest number
$d$ such that $T$ can be written as a sum of $d$ tensors of the form $v^{\otimes d}=v\otimes v\otimes \cdots \otimes v$. It is clear that the tensor rank is less than or equal to the symmetric rank. It was conjectured by Comon \cite{Comon} that the symmetric rank and tensor rank of a symmetric tensor are equal. Recently, Shitov gave a counterexample \cite{Shitov}. In this paper, we study the notion of $G$-stable rank of a tensor. The $G$-stable rank of a tensor is defined in terms of geometric invariant theory
and the notion of stability for algebraic group actions on tensors.
It is also natural to define a symmetric $G$-stable rank for a symmetric tensor. One main result of this paper is that the symmetric $G$-stable rank and $G$-stable rank of a symmetric tensor are the same. 

In algebraic geometry and singularity theory, the log canonical threshold is an important invariant of singularities. We will show that the symmetric $G$-stable rank and the log canonical threshold are closely related. We extend the notion of $G$-stable rank to ideals in a coordinate ring of a smooth irreducible affine variety. In this context, we show 
 that the log canonical threshold is less than or equal to the $G$-stable rank. In the case of monomial ideals in the polynomial ring we show equality.

\subsection{Stability of tensors}
Let $K$ be a perfect field, and $V$ a finite dimensional vector space over $K$, we consider the action of the group of product of special linear groups $\SL(V)^d = \SL(V)\times \SL(V) \times \cdots \times \SL(V)$ on the tensor product space $V^{\otimes d}=V\otimes V\otimes\cdots\otimes V$. A 1-parameter subgroup of an algebraic group $G$ is a homomorphism of algebraic groups $\lambda:{\mathbb G}_m\to G$, where $\mathbb G_m$ is the multiplicative group. For any integer $m$, we define the multiple of $\lambda$ by $m$, denoted by $m\cdot\lambda$, which is also a 1-parameter subgroup with $(m\cdot\lambda)(t) = (\lambda(t))^m$. A 1-parameter subgroup is {\em indivisible} if it is not a multiple of any other 1-parameter subgroup with factor $m\geq 2$. We say a tensor $v\in V^{\otimes d}$ is $\SL(V)^d$-{\em unstable} if there is a 1-parameter subgroup $\lambda: \mathbb G_m\to \SL(V)^d$, such that 
\[
\lim_{t\to 0}\lambda(t)\cdot v = 0.
\]

If no such 1-parameter subgroup exists, then $v$ is called $\SL(V)^d$-{\em semistable}. Let $G$ be a reductive algebraic group over $K$. By a $G$-scheme $X$ we mean a separated, finite type scheme $X$ over $K$ as well as a morphism $G\times X\to X$ mapping $(g,x)$ to $g\cdot x$, such that $g\cdot(h\cdot x) = (gh)\cdot x$, for all $g,h\in G$ and for all $x\in X$. A morphism $f: X\to Y$ between two $G$-schemes $X$ and $Y$ is $G$-equivariant if  for all $g\in G$
and $x\in X$, we have $f(g\cdot x) = g\cdot f(x)$. A subscheme $S$ of a $G$-scheme $X$ is called a $G$-subscheme if $S$ is a $G$-scheme and the immersion $S\hookrightarrow X$ is $G$-equivariant. 

Throughout this paper, we will work over a perfect field $K$. In \cite{Kempf}, Kempf proved a $K$-rational version of the \textit{Hilbert-Mumford} criterion:

\begin{theorem}[\cite{Kempf}, Corollary 4.3]
\label{HM}
 Let $G$ be a reductive algebraic group. Suppose that  $X$ is a $G$-scheme and $x\in X$ is a $K$-point. Assume $S$ is a closed $G$-subscheme of $X$ which does not contain $x$ and $S$ meets the closure of the orbit $G\cdot x$. Then there exits a $K$-rational 1-parameter subgroup $\lambda :\mathbb G_m\to G$, such that
\[
\lim_{t\to 0}\lambda(t)\cdot x\in S.
\]
\end{theorem}

If $G = \SL(V)^d$, $X = V^{\otimes d}$ , then by Theorem \ref{HM}, $v$ is unstable if and only if $0$ is in the closure of the orbit of $v$, i.e. $0\in \overline{\SL(V)^d\cdot v}$. 

A tensor $T\in V^{\otimes d}$ is called symmetric
if it is invariant under the action of symmetric group $S_d$.
Let $D^d V\subseteq V^{\otimes d}$ be the space
of symmetric tensors. As a representation of $\GL(V)$, 
this is the space of divided powers,
which is isomorphic to the $d$-th symmetric power $S^dV$ if the characteristic of $K$ is $0$ or $>d$.

It is  interesting to look at the diagonal action of $\SL(V)$  on $\symm^dV$ via the diagonal embedding: 
\begin{equation}
\Delta: \SL(V)\hookrightarrow \SL(V)^d.
\end{equation}

\begin{definition}\label{ss}
A symmetric tensor $v\in \symm^dV$ is $\SL(V)$-\textit{unstable} if there is a 1-parameter subgroup $\lambda: \mathbb G_m\to \SL(V)$, such that 
\[
\lim_{t\to 0}\lambda(t)\cdot v = 0.
\]
Otherwise we say $v$ is $\SL(V)$-{\em semistable}.
\end{definition}

 
 

 \subsection{G-stable rank for tensors}

In \cite{G-stable rank}, Derksen introduced $G$-stable rank for tensors. Suppose the base field $K$ is perfect. If $\lambda:\Gm\to \GL_n$ is a $1$-parameter subgroup,
then we can view $\lambda(t)$ as an invertible $n\times n$ matrix whose entries lie in the ring $K[t,t^{-1}]$ of Laurent polynomials. We say that $\lambda(t)$ is a polynomial $1$-parameter subgroup of $\GL_n$ if all these entries lie in the polynomial ring $K[t]$. 
 Consider the action of the group $G = \GL(V_1)\times \GL(V_2)\times\cdots\times \GL(V_d)$ on the tensor product space $W = V_1\otimes V_2\otimes\cdots\otimes V_d$. 
A 1-parameter subgroup $\lambda:\Gm\to G$ 
can be written as
\[
\lambda(t) = (\lambda_1(t), \cdots, \lambda_d(t)),
\]
where $\lambda_i(t)$ is a $1$-parameter subgroup of $\GL(V_i)$ for all $i$. We say that $\lambda(t)$ is polynomial if and only if $\lambda_i(t)$ is a polynomial $1$-parameter subgroup for all $i$.

The $t$-valuation $\val(a(t))$ of a polynomial $a(t)\in K[t]$ is the biggest integer $n$ such that $a(t) = t^nb(n)$ for some $b(t)\in K[t]$. For $a(t), b(t)\in K[t]$, the $t$-valuation $\val\big(\frac{a(t)}{b(t)}\big)$ of the rational function $\frac{a(t)}{b(t)}\in K(t)$ is $\val\big(\frac{a(t)}{b(t)}\big) = \val(a(t)) - \val(b(t))$. For a tuple $u(t) = (u_1(t),u_2(t),\cdots, u_d(t))\in K(t)^d$, we define the $t$-valuation of $u(t)$ as
\begin{equation}\label{fval}
\val(u(t)) = \min_i\{\val(u_i(t))| 1\le i\le d\}.
\end{equation}

If $\lambda$ is a $1$-parameter subgroup of $G$ and $v\in W$ is a tensor, then we have $\lambda(t)\cdot v\in K(t)\otimes W$. We view $K(t)\otimes W$ as a vector space over $K(t)$ and define the $t$-valuation $\val(\lambda(t)\cdot v)$ as in \eqref{fval}. Assume $\val(\lambda(t)\cdot v) > 0$, then for for any $\alpha = (\alpha_1,\alpha_2,\cdots,\alpha_d)\in \mathbb R^d_{>0}$, we define the slope
\begin{equation}
    \mu_\alpha(\lambda(t), v) = \frac{\sum_{i = 1}^d\alpha_i \val(\det(\lambda_i(t)))}{\val(\lambda(t)\cdot v)}.
\end{equation}

The $G$-stable rank for $v\in W$ is the infimum of the slope with respect to all such 1-parameter subgroups. More precisely:

\begin{definition}[\cite{G-stable rank}, Theorem 2.4]
If $\alpha\in \mathbb R^d_{>0}$, then the $G$-stable rank $\rk^G_\alpha(v)$ is the infimum of $\mu_\alpha(\lambda(t), v)$ where $\lambda(t)$ is a polynomial 1-parameter subgroup of $G = \GL(V_1)\times\cdots\times \GL(V_n)$ and $\val(\lambda(t)\cdot v) > 0$. If $\alpha = (1,1,\cdots, 1)$, we simply write $\rk^G(v)$.
\end{definition}

Let $V$ be a finitely dimensional vector space over $K$. Let $\symm^dV\subset V^{\otimes d}$ be the space of all symmetric tensors. Assume the group $\GL(V)$ acts on $\symm^dV$ via the diagonal embedding: $\GL(V) \ \hookrightarrow \GL(V)^d$.
\begin{definition}
Let $v\in \symm^dV$ be a symmetric tensor, the symmetric $G$-stable rank $\symmrk^G(v)$ of $v$ is the infimum of $\mu(\lambda(t),v) = d\frac{\val(\det(\lambda(t)))}{\val(\lambda(t)\cdot v)}$, where $\lambda(t)$ is a polynomial 1-parameter subgroup of $\GL(V)$ and $\val(\lambda(t)\cdot v) > 0$.
\end{definition}
Since any 1-parameter subgroup of $\GL(V)$ is also a 1-parameter subgroup of $\GL(V)^d$ via the diagonal embedding, we have $\symmrk^G(v)\ge \rk^G(v)$ for any $v\in \symm^d V$.  It turns out that the other inequality is also true,
\begin{theorem}\label{second result}
Let $v\in \symm^d V$ be a symmetric tensor, then we have
\[
\symmrk^G(v) = \rk^G(v).
\]
\end{theorem}

\begin{example}
Suppose that $V=K^2$, and $v = e_2\otimes e_1\otimes e_1+e_1\otimes e_2\otimes e_1+e_1\otimes e_1\otimes e_2\in V^{\otimes 3}$, where $\{e_1,e_2\}$ is the standard basis of $V=K^2$. Let $\lambda(t) = \begin{pmatrix}
t&0\\
0 &1
\end{pmatrix}$ be a polynomial 1-parameter subgroup of $\GL(K^2)$. 
Then $\lambda(t)\cdot v = t^2 v$, $\det(\lambda(t)) = t$, the slope is
\[
\mu(\lambda(t), v) = 3\frac{\val(\det(\lambda(t)))}{\val(\lambda(t)\cdot v)}=\frac{3}{2}.
\]

Therefore we have $\symmrk^G(v)\le \frac{3}{2}$. It was proved in \cite{G-stable rank} that $\rk^G(v)=\frac{3}{2}$. Hence by the fact $\symmrk^G(v)\ge \rk^G(v)$ we have $\symmrk^G(v) = \frac{3}{2}$.

\end{example}

A 1-parameter subgroup of $\SL(V)$ is also a 1-parameter subgroup of $\SL(V)^d$ via the diagonal embedding. It follows that if a symmetric tensor $v\in \symm^d V$ is $\SL(V)$-unstable, 
then $v$ is also $\SL(V)^d$-unstable. Equivalently, if $v$ is $\SL(V)^d$-semistable, then $v$ is also $\SL(V)$-semistable. It follows from Theorem \ref{second result} that the converse direction is also true:

 \begin{corollary}\label{first result}
 Let $v\in \symm^d V$ be a symmetric tensor, then $v$ is $\SL(V)^d$-semistable if and only if it is $\SL(V)$-semistable.
 \end{corollary}



\subsection{G-stable rank for ideals and log canonical threshold}

Let $V$ be an $n$-dimensional vector space over a perfect field $K$. By choosing a basis of $V$ and a dual basis $\{x_1,x_2,\dots,x_n\}$ of $V^\star$, we have an isomorphism of algebras $SV^\star\cong K[x_1,\cdots, x_n]$, where $SV^\star$ is the symmetric algebra on the vector space $V^\star$. 
We have defined the symmetric $G$-stable rank for symmetric tensors, it is natural to extend this idea to polynomials and more generally to ideals in the polynomial ring $K[x_1,\cdots,x_n]$. Furthermore, let $X$ be a smooth irreducible affine variety with coordinate ring $K[X]$, and let $\mathfrak a\subset K[X]$ be an ideal. We can define the $G$-stable rank $\rk^G(P, \mathfrak a) $ for the ideal $\mathfrak a$ at a point $P\in V(\mathfrak a)$. We postpone the precise definition of $G$-stable rank for ideals to Section~\ref{ideal rank}. It turns out that the $G$-stable rank $\rk^G(P, \mathfrak a)$ of an ideal $\mathfrak a$ at $P$ is closely related to the {\em log canonical threshold} $\lct_P(\mathfrak a)$ of the ideal $\mathfrak a$ at the point $P\in V(\mathfrak a)$.

Log canonical threshold is an invariant of singularities in algebraic geometry, \cite{Mus} gives a comprehensive introduction to this subject.  Let $K = \mathbb C$ be the complex field. Let $H\subset \mathbb C^n$ be a hypersurface defined by a polynomial $f\in \mathbb C[x_1,\cdots,x_n]$, and let $P\in H$ be a closed point. The log canonical threshold $\lct_P(f)$ of $f$ at the point $P$ tells us how singular $f$ is at the point $P$. More precisely, $\lct_P(f)$ is a rational number bounded above by 1, and equal to 1 if $P$ is a smooth point of $H$.

There are several equivalent ways to define the log canonical threshold, here we give an analytic definition, which we will use later. 
\begin{definition}\label{lct}
Let $X$ be a smooth irreducible affine variety. Let $\mathfrak a = (f_1,\cdots,f_r) \subset \mathbb C[X]$ be an ideal, and $P\in V(\mathfrak a)$ is a closed point.
The {\em log canonical threshold} $\lct_P(\mathfrak a)$ of the ideal $\mathfrak a$ at $P$ is 
\begin{equation}
\lct_P(\mathfrak a) = \sup\Big\{s>0\ \Big|\ \frac{1}{(\sum_{i = 0}^r|f_i|^2)^s}  \
    \text{is}\ \text{integrable}\ \text{around}\ P\Big\}.
\end{equation}
\end{definition}

The log canonical threshold $\lct_P(f)$ of a polynomial $f\in \mathbb C[X]$ is the log canonical threshold of the principle ideal $\mathfrak a = (f)$. We have the following relation between the log canonical threshold and the $G$-stable rank:


\begin{theorem}\label{relation}
In the situation of Definition \ref{lct}, the log canonical threshold of $\mathfrak a$ is less than or equal to the $G$-stable rank of $\mathfrak a$ at $P$:
\begin{equation}
\lct_P(\mathfrak a) \le \rk^G(P,\mathfrak a).
\end{equation}
\end{theorem}

When $\mathfrak a$ is a monomial ideal, i.e. $\mathfrak a$ is generated by monomials, the equality holds.

\begin{theorem}\label{equal}
Suppose $\mathfrak a\subset \mathbb C[x_1,\cdots,x_n]$ is a proper nonzero ideal generated by monomials and $P=(0,\cdots,0)$ is the origin. Then we have
\begin{equation}
\lct_P(\mathfrak a) = \rk^G(P, \mathfrak a).
\end{equation}
\end{theorem}


\section{ Kempf's theory of optimal subgroups}\label{Kempf Results}
Let $G$ be a reductive algebraic group over a perfect field $K$. In our case, $G$ is one of $\SL(V)^d, \SL(V), \GL(V)^d$ or $\GL(V)$, depend on the situation. Let $\Gamma(G)$ denote the set of all 1-parameter subgroups of $G$.  
In \cite{Kempf}, Kempf provided a way to approach the boundary of an orbit. Following~\cite{Kempf}, we have the definition:

\begin{definition}
Let $X$ be a $G$-scheme over a perfect field $K$ and let $x\in X$ be a $K$-point. We define $ |X,x|$ to be the set of all 1-parameter subgroups of $G$ such that $\lim_{t\to 0}\lambda(t)\cdot x$ exists in $X$. 
Assume $S$ is a $G$-invariant closed sub-scheme of X not containing $x$, we define a subset $|X,x|_S\subset |X,x|$ by 
\begin{equation}
|X,x|_S =\{\lambda\in |X,x|\mid\lim_{t\to 0}\lambda(t)\cdot x \in S\}. 
\end{equation}

\end{definition}

\begin{remark}
If $S\cap \overline{G\cdot x}\ne \emptyset$, then by Theorem \ref{HM}, there exists a 1-parameter subgroup $\lambda(t)\in \Gamma(G)$, such that $\lim_{t\to 0}\lambda(t)\cdot x\in S$. Hence $|X,x|_S\ne \emptyset$. If $X = V^{\otimes d}$, $S = 0$ and $v\in V^{\otimes d}$ is $\SL(V)^d$-unstable, then $|V^{\otimes d}, v|_{\{0\}}\ne \emptyset$.
\end{remark}

Let $\lambda\in |X,x|$ be a 1-parameter subgroup of $G$, we get a morphism $\phi_\lambda:\mathbb A^1\to X$ by $\phi_\lambda(t) = \lambda(t)\cdot x$ if $t\neq 0$ and $\phi_\lambda(0)=\lim_{t\to 0} \lambda(t)\cdot x$. Assume $S$ is a $G$-invariant closed sub-scheme of X not containing $x$, the inverse image $\phi_\lambda^{-1}(S)$ is an effective divisor supported inside $t=0$. Let $a_{S,x}(\lambda)$ denote the degree of the divisor $\phi_\lambda^{-1}(S)$ for $\lambda\in |X,x|$. Note that we have a natural conjugate action of $G$ on the set of 1-parameter subgroups $\Gamma(G)$ by $(g\cdot \lambda)(t) = g\lambda(t)g^{-1}$, where $g\in G$, $\lambda\in \Gamma(G)$.

\begin{definition}\label{def:length}
A {\it length} function $\|\cdot\|$ is a non-negative real-valued function on $\Gamma(G)$ such that
\begin{enumerate}
    \item $\|g\cdot\lambda\| = \|\lambda\|$ for any $\lambda\in\Gamma(G)$ and $g\in G$.
    \item For any maximal torus $T \subseteq G$, we have $\Gamma(T)\subseteq \Gamma(G)$, the restriction of $\|\cdot\|$ on $\Gamma(T)$ is integral valued and extends to a {\it norm}  on the vector space $\Gamma(T)\otimes_\mathbb Z \mathbb R$.
  \end{enumerate}
\end{definition} 

\begin{remark}\label{len_exists}
Such a length function exists. Let $T$ be a maximal torus of $G$. Let $N$ be the normalizer of $T$. Then the Weyl group with respect to $T$ is defined by $W = N/T$. By the fact that $\Gamma(G)/G\cong \Gamma(T)/W$, it suffices to define a $W$-invariant norm on $\Gamma(T)\otimes_\mathbb Z \mathbb R$. Since $W$ is a finite group, any norm on $\Gamma(T)\otimes_\mathbb Z \mathbb R$ and then average over $W$ will work.
\end{remark}

\begin{remark}
In the original paper \cite{Kempf}, Kempf defined a length function $\|\cdot\|$ that satisfies a different condition (2): for any maximal torus $T$ of $G$, there is a positive definite integral-valued bilinear form $(\ ,\ )$ on $\Gamma(T)$, such that $(\lambda,\lambda) = \|\lambda\|^2$ for any $\lambda$ in $\Gamma(T)$. But the proof in \cite{Kempf} of the theorem below 
is also valid for our slightly weaker definition of length function.

\end{remark}


\begin{theorem}[Kempf \cite{Kempf}]\label{main theorem}
Let $X$ be an affine $G$-scheme over a perfect field $K$. Let $x\in X$ be a $K$-point. Assume $S$ is an $G$-invariant closed sub-scheme not containing $x$ such that $S\cap \overline{G\cdot x} \ne \emptyset$. Fix a length function $\|\cdot\|$ on $\Gamma(G)$, then we have
\begin{enumerate}
    \item The function $\frac{a_{S,x}(\lambda)}{\|\lambda\|}$ has a maximum positive value $B_{S,x}$ on the set of non-trivial 1-parameter subgroups in $|X,x|$.
    \item Let $\Lambda_{S,x}$ be the set of indivisible 1-parameter subgroups $\lambda\in |X,x|$ such that $a_{S,x}(\lambda)=B_{S,x}\cdot \|\lambda\|$, then we have
    \begin{enumerate}
        \item $\Lambda_{S,x}\ne\emptyset$.
        \item For $\lambda\in \Lambda_{S,x}$, Let $P(\lambda)=\{g\in G|\lim_{t\to 0} \lambda(t)\cdot g\cdot\lambda(t)^{-1} \  \mbox{exists}\}$, then $P(\lambda)$ is a parabolic subgroup and independent of $\lambda$. We denote it by $P_{S,x}$.
        \item Any maximal torus of $P_{S,x}$ contains a unique member of $\Lambda_{S,x}$.
    \end{enumerate}
\end{enumerate}
\end{theorem}

\section{$G$-stable rank and symmetric $G$-stable rank}\label{equal ranks}
Let $V$ be an $n$ dimensional vector space over $K$, fix a maximal torus $T$ of $\GL(V)$, we have an isomorphism $\Gamma(T)\cong \mathbb Z^n$. Any 1-parameter subgroup $\lambda$ of the maximal torus $T$ is given by a tuple of $n$ integers $(\nu_1,\cdots,\nu_n)$, we define a function on $\Gamma(T) \cong \mathbb Z^n$ by
\begin{equation}\label{len1}
\|\lambda\| = \sum_{i = 1}^n|\nu_i|.
\end{equation}
This function extends linearly to a norm on the vector space $\Gamma(T)\otimes_{\mathbb Z} \mathbb R$.
The Weyl group of $\GL(V)$ with respect to $T$ is the symmetric group $S_n$. It is clear that the function is invariant under the action of $S_n$ by permutation, therefore by Remark \ref{len_exists}, it defines a length function on $\Gamma(\GL(V))$. Let $G = \GL(V)^d$, fix a maximal torus $T_i\subset\GL(V)$ for each component of $\GL(V)^d$, then $T = T_1\times\cdots\times T_d$ is a maximal torus of $G$. We have $\Gamma(T) \cong (\mathbb Z^n)^d$, Let $\lambda = (\lambda_1,\cdots,\lambda_d)$ be a 1-parameter subgroup of $T$, where 
\[
\lambda_i=(\lambda_{i,1},\cdots,\lambda_{i, n}), \ \lambda_{i, j} \in \Z\ \text{for\ all}\  j
\]
is a tuple of $n$ integers. 
 We define a function on $\Gamma(T)$ by 
 
 \begin{equation}\label{len2}
 \|\lambda\| = \sum_{i = 1}^d \|\lambda_i\|,
 \end{equation}
 where $\|\lambda_i\| = \sum_{j = 1}^n |\lambda_{i, j}|$. This extends to a length function on $\Gamma(G) = \Gamma(\GL(V)^d)$.
 
 Let $G=\GL(V)^d$, $X= V^{\otimes d}$ and $S=\{0\}$. Recall the definition of $t$-valuation in equation~(\ref{fval}).

\begin{lemma}
Let $v\in \symm^dV\subset V^{\otimes d}$ be a symmetric tensor, and $G = \GL(V)^d$ acts on $V^{\otimes d} $ in the usual way. If $\lambda(t)$ is a 1-parameter subgroup of $G$, then
\begin{enumerate}
    \item $|X,v| = \{\lambda\in \Gamma(G)\mid \val(\lambda(t)\cdot v)\ge 0\}$.
    \item $|X,v|_{\{0\}} = \{\lambda\in \Gamma(G)\mid \val(\lambda(t)\cdot v)> 0\}$.
    \item $a_{\{0\},v}(\lambda) = \val(\lambda(t)\cdot v)$ for $\lambda\in |X,v|$.
\end{enumerate}
\end{lemma}
 
 \begin{proof}
 This follows immediately from the definition.
 \end{proof}
\begin{lemma}\label{>=0}
Let $v\in \symm^dV\subset V^{\otimes d}$ be a symmetric tensor, then the function $\frac{\val(\lambda(t)\cdot v)}{\|\lambda\|}: \Gamma(G)\to \R $ attains its maximal value at some 1-parameter subgroup $\lambda\in \Gamma(\GL(V)^d)$. There exists a maximal torus $T\subset \GL(V)^d$, such that $\lambda\in \Gamma(T)$ and under the isomorphism $\Gamma(T^d)\cong(\mathbb Z^n)^d$, we can write $\lambda = (\lambda_1,\cdots,\lambda_d)$, where $\lambda_i =(\lambda_{i,1},\cdots,\lambda_{i,n})$ such that $\lambda_{i,j}\in \Z$ and $\lambda_{i, j} \ge 0$ for all $1\le i\le d, 1\le j\le n$, in other words, $\lambda$ is a polynomial 1-parameter subgroup.
\end{lemma}

\begin{proof}

By Theorem \ref{main theorem}, the maximal value of $\frac{\val(\lambda(t)\cdot v)}{\|\lambda\|}$ exists. Let $T\subset G$ be a maximal torus and $\lambda\in \Gamma(T)$ such that the function attains its maximum at $\lambda$. Assume $\lambda$ is of the form in the lemma and $\lambda_{i, j} < 0$ for some $i$ and $j$. If we replace $\lambda_{i, j}$ by $-\lambda_{i, j}$, $\val(\lambda(t)\cdot v)$ never decreases and $\|\lambda\|$ does not change. Therefore by the maximality of $\frac{\val(\lambda(t)\cdot v
)}{\|\lambda\|}$, the value $\frac{\val(\lambda(t)\cdot v)}{\|\lambda\|}$ does not change after the replacement. Hence without loss of generality we can assume all $\lambda_{i, j} \ge 0$.
\end{proof}

Recall that for a tensor $v\in V^{\otimes d}$ and a polynomial 1-parameter subgroup $\lambda$ of $G = \GL(V)^d$ such that $\val(\lambda(t)\cdot v)>0$, we have the slope function
\begin{equation}\label{slope}
\mu(\lambda(t), v) = \frac{\sum_{i = 1}^d \val(\det(\lambda_i(t)))}{\val(\lambda(t)\cdot v)}.
\end{equation}

Let $\lambda = (\lambda_1,\cdots, \lambda_d)\in \Gamma(G)$ be a polynomial 1-parameter subgroup of $G=\GL(V)^d$, then by Lemma \ref{>=0}, $\sum_{i=1}^d\val(\det(\lambda_i(t)))$ is the restriction of the length function defined by equation (\ref{len2}). Let $S_d$ be the symmetric group acting on $G=\GL(V)^d$ by permuting the $d$ components. Then the length function defined by equation (\ref{len2}) is invariant under the action of $S_d$. From now on, fix this length function on $\Gamma(G)$. We have a corollary following from Theorem \ref{main theorem}:

\begin{corollary}\label{Pd}
Let $v\in \symm^d V\subset V^{\otimes d}$ be a symmetric tensor. Let $\Lambda_{\{0\},v}$ be the set of indivisible 1-parameter subgroups $\lambda\in |V^{\otimes d},v|$ such that $\frac{\val(\lambda(t)\cdot v)}{\|\lambda\|}$ attains the maximum value. Then we have
\begin{enumerate}
    \item $\Lambda_{\{0\},v}$ is invariant under $S_d$.
    \item $P_{\{0\},v}$ is $S_d$ invariant. In other words, $P_{\{0\},v} = P^d\subset \GL(V)^d$ for some parabolic subgroup $P\subset\GL(V)$.
    
\end{enumerate}
\end{corollary}

\begin{proof}\ 
\begin{enumerate}
    \item It is clear that $\frac{\val(\lambda(t)\cdot v)}{\|\lambda\|}$ is $S_d$ invariant. Indeed, Let $\sigma\in S_d$, since $v\in \symm^d V$ is a symmetric tensor and $\|\cdot\|$ is $S_d$ invariant, we have
    
    \[
    \frac{\val((\sigma \lambda)(t)\cdot v)}{\|\sigma \lambda\|} = \frac{\val((\sigma \lambda)(t)\cdot (\sigma v))}{\|\sigma \lambda\|} =
    \frac{\val(\sigma (\lambda(t)\cdot v))}{\|\sigma \lambda\|} = \frac{\val(\lambda(t)\cdot v)}{\|\lambda\|}.
    \]
    Therefore if $\lambda\in \Lambda_{\{0\},v}$, so is $\sigma(\lambda)$.

    \item  
    Since $G = \GL(V)^d$, the parabolic subgroup $P_{\{0\},v}$ is a product of parabolic subgroups of $\GL(V)$, the symmetric group $S_d$ acts on $P_{\{0\},v}$ by permuting the components. Let $\lambda\in \Lambda_{\{0\},v}$, for any $\sigma\in S_d$, we have
    \[
    \sigma(P_{\{0\},v}) = P(\sigma(\lambda)) = P_{\{0\},v}.
    \]
    We used the fact that $P_{\{0\},v}=P(\lambda)$ is independent of $\lambda\in\Lambda_{\{0\},v}$ and $\sigma(\lambda)\in \Lambda_{\{0\},v}$. So $P_{\{0\},v}$ is $S_d$ invariant and we can find a parabolic subgroup $P\subset \GL(V)$ such that $P_{\{0\},v} = P^d$.
\end{enumerate}
\end{proof}

Let $T\subset P\subset \GL(V)$ be a maximal torus, then $T^d$ is a maximal torus of $P^d=P_{\{0\}, v}$. By (2.c) in Theorem \ref{main theorem} and Lemma \ref{>=0}, there is a polynomial 1-parameter subgroup $\lambda = (\lambda_1,\cdots, \lambda_n)$ of $T^d\subset \GL(V)^d$, such that the slope function 
\[
\mu(\lambda(t),v)=\frac{\sum_{i = 1}^d \val(\det(\lambda_i(t)))}{\val(\lambda(t)\cdot v)} = \frac{\|\lambda\|}{\val(\lambda(t)\cdot v)}
\]
has a minimum value at $\lambda$. The minimal value of $\mu(\lambda(t), v)$ is by definition the $G$-stable rank $\rk^G(v)$ of $v$. In rest of the section, we fix such a maximal torus $T\subset \GL(V)$.

 Let $\lambda=(\lambda_1,\cdots,\lambda_d)$ be a polynomial 1-parameter subgroup of $T^d\subset \GL(V)^d$, then $\gamma = \Pi_{i=1}^d\lambda_i$ is a polynomial 1-parameter subgroup of $\GL(V)$. Furthermore, $\gamma$ acts on $v\in \symm^dV$ via the diagonal embedding $\GL(V)\hookrightarrow \GL(V)^d$. We have the following lemma:

\begin{lemma}
For any symmetric tensor $v\in \symm^d V$, we have $\val(\gamma(t)\cdot v) \ge d\cdot \val(\lambda(t)\cdot v)$.
\end{lemma}

\begin{proof}
Let $C = \val(\lambda(t)\cdot v)$, we define a subspace $W$ of $V^{\otimes}$ as following
\[
W = \{w\in V^{\otimes d}|\val(\sigma(\lambda(t))\cdot w)\ge C, \forall \sigma\in S_d\}.
\]
Since $\val(\sigma(\lambda(t))\cdot v) = \val(\sigma(\lambda(t)\cdot v)) = \val(\lambda(t)\cdot v) = C$, we have $v\in W$. For any $\sigma\in S_d$, it is clear that $\sigma(\lambda(t))\cdot W\subset t^CK[t]\cdot W$. We can write
\begin{align*}
\gamma(t)\cdot v
&=(\Pi_{i=1}^d\lambda_i,\cdots,\Pi_{i=1}^d\lambda_i)\cdot v
\\
&=(\lambda_1,\lambda_2,\cdots,\lambda_d)(\lambda_2,\lambda_3,\cdots,\lambda_d,\lambda_1)\cdots(\lambda_d,\lambda_1,\cdots,\lambda_{d-1})\cdot v
\\
&=(\lambda_1,\lambda_2,\cdots,\lambda_d)\sigma(\lambda_1,\lambda_2,\cdots,\cdots,\lambda_d)\cdots\sigma^{d- 1}(\lambda_1,\lambda_2,\cdots,\lambda_{d})\cdot v
\\
&=\lambda\sigma(\lambda)\cdots\sigma^{d-1}(\lambda)\cdot v,
\end{align*}
where $\sigma\in S_d$ satisfies $\sigma(1) = 2, \sigma(2) = 3,\cdots,\sigma(d) = 1$. Therefore $\gamma(t)\cdot v \in t^{dC}K[t]\cdot W$, hence we get $\val(\gamma(t)\cdot v) \ge dC = d\cdot \val(\lambda(t)\cdot v)$.
\end{proof}

Next we prove that the symmetric $G$-stable rank is the same as the $G$-stable rank for symmetric tensors.


\begin{proof}[Proof of Theorem \ref{second result}]
Let $T$ be the chosen maximal torus of $\GL(V)$ as above. Let $\lambda = (\lambda_1,\cdots,\lambda_d)$ be a polynomial 1-parameter subgroup of $T^d\subset \GL(V)^d$ such that the slope function $\mu(\lambda(t),v)$ attains its minimum value. In other words, $\lambda$ computes the $G$-stable rank $\rk^G(v)$ of $v$:
\[
\rk^G(v) = \frac{\sum_{i=1}^d \val(\det(\lambda_i(t)))}{\val(\lambda(t)\cdot v)}.
\]
Let $\gamma = \Pi_{i=1}^d\lambda_i$ as above, then we have
\[
\symmrk^G(v)\le \frac{\sum_{i=1}^d \val(\det(\gamma(t)))}{\val(\gamma(t)\cdot v)} \le  \frac{d\sum_{i=1}^d \val(\det(\lambda_i(t)))}{d\cdot \val(\lambda(t)\cdot v)} = \rk^G(v).
\]
On the other hand, it is clear that $\symmrk^G(v)\ge \rk^G(v)$. Therefore $\symmrk^G(v) = \rk^G(v)$, this completes the proof.
\end{proof}

\section{Stability of symmetric tensors}\label{stability}

As a result of Theorem \ref{second result}, we prove Corollary \ref{first result}, which says that for a symmetric tensor $v\in \symm^d V$, $v$ is $\SL(V)^d$-semistable if and only if $v$ is $\SL(V)$-semistable. It is clear that $\SL(V)^d$-semistability implies $\SL(V)$-semistability. To prove the other direction, we will use a result which relates semistability with $G$-stable rank.

\begin{prop}[\cite{G-stable rank}, Proposition 2.6]
\label{2.6}
Suppose that $\alpha=(\frac{1}{n_1},\cdots\frac{1}{n_d})$ where $n_i =\dim V_i$. For $v\in V_1\otimes V_2\otimes \cdots\otimes V_d$ we have $\rk^G_\alpha(v) \le 1$. Moreover, $\rk^G_\alpha(v) = 1$ if and only if $v$ is semistable with respect to the group $H = \SL(V_1)\times\SL(V_2)\times\cdots\times\SL(V_d)$.
\end{prop}

If $v\in \symm^d V$ is a symmetric tensor, $\alpha = (1,1,\cdots,1)$ and $n = \dim V$, then by the above proposition, $\rk^G(v) = n$ if and only if $v$ is $\SL(V)^d$-semistable. We have a similar result for symmetric $G$-stable rank.

\begin{prop}\label{2.66}
For a symmetric tensor $v\in \symm^d V$, we have $\symmrk^G(v) \le n$, where  $n = \dim V$. Moreover, $\symmrk^G(v) = n$ if and only if $v$ is $\SL(V)$-semistable.
\end{prop}
\begin{proof}
The first statement is clear from Theorem \ref{second result}. If $\symmrk^G(v) = n$, by Proposition \ref{2.6}, we have $\rk^G(v) = n$ and $v$ is $\SL(V)^d$-semistable, hence $v$ is $\SL(V)$-semistable. On the other hand,  assume $v$ is $\SL(V)$-semistable. Let $\lambda$ be a polynomial 1-parameter subgroup of $\GL(V)$ such that $\lim_{t\to 0}\lambda(t)\cdot v=0$. Then we can define another 1-parameter subgroup $\lambda'(t) = \lambda(t)^n t^{-e}$, where $\det(\lambda(t)) = t^e$, such that $\det(\lambda')=1$ and $\lambda'\in \SL(V)$. Since $v$ is $\SL(V)$-semistable, we have $\val(\lambda'(t)\cdot v) \le 0$. It follows that
\[
\val(\lambda'(t)\cdot v) = \val(\lambda(t)^d t^{-e}\cdot v) = n\val(\lambda(t)\cdot v) - ed\le 0.
\]
The slope function
\[
\mu(\lambda(t), v) = \frac{d\val(\det(\lambda(t)))}{\val(\lambda(t)\cdot v)} = \frac{de}{\val(\lambda(t)\cdot v)}\ge n.
\]
We get $\symmrk^G(v) = n$.
\end{proof}

\begin{proof}[Proof of Corollary \ref{first result}]
It suffices to prove that if $v$ is $\SL(V)$-semistable, then $v$ is $\SL(V)^d$-semistable. Let us assume $v$ is $\SL(V)$-semistable, then by Proposition \ref{2.66} and Theorem \ref{second result}, 
\[
\rk^G(v) = \symmrk^G(v) = n.
\]
It follows from Proposition \ref{2.6} that $v$ is $\SL(V)^d$-semistable.
\end{proof}

\section{G-stable rank for ideals and log canonical threshold} \label{ideal rank}

\subsection{G-stable rank for ideals}
Let $X = \text{Spec}(R)$ be a nonsingular irreducible complex affine algebraic variety of dimension $n$, and $\mathfrak a\subset R$ be a nonzero ideal, and let $P\in V(\mathfrak a)$ be a closed point, $\mathcal O_P$ be the local ring at $P$ and $\mathfrak m_P$ be the maximal ideal corresponding to $P$. 

\begin{definition}[\cite{Sha}]
Functions $x_1,\cdots,x_n\in \mathcal O_P$ are {\em a system of local parameters} at $P$ if each $x_i\in \mathfrak m_P$, and the images of $x_1,\cdots,x_n$ form a basis of the vector space $\mathfrak m_P/\mathfrak m_P^2$.
\end{definition}

Let $T = \{x_1, x_2,\cdots, x_n\}$ be a system of local parameters at $P$.  Let $\C\{x_1,x_2,\dots,x_n\}$ be the ring of convergent power series in $x_1,x_2,\dots,x_n$. The ring $\mathcal O_P$
is contained in $\C\{x_1,x_2,\dots,x_n\}$.
If $y_1,y_2,\dots,y_n$ is any system of local parameters, then $\C\{x_1,x_2,\dots,x_n\}=\C\{y_1,y_2,\dots,y_n\}$. For any $\lambda = (\lambda_1,\lambda_2,\cdots,\lambda_n)\in \mathbb Z^n_{\ge 0}$, we have a natural action of $\mathbb C^*$ on $\C\{x_1,x_2,\cdots,x_n\}$ by $t\cdot x_i = t^{\lambda_i}x_i$ for any $t\in\mathbb C^*$.

\begin{definition}\label{val&ord}
Let $T = \{x_1,x_2,\cdots,x_n\}$ be a system of local parameters at $P$ and $\lambda = (\lambda_1,\lambda_2,\cdots,\lambda_n)\in \mathbb Z^n_{\ge 0}$ a tuple of non-negative integers. Let $f\in \C\{x_1,x_2,\cdots,x_n\}$ be a convergent power series. We define the {\em valuation} of $f$ with respect to $T$ by
\begin{equation}\label{val}
    \val^T_\lambda(f) = \max\{k|f(t^{\lambda_1}x_1,t^{\lambda_2}x_2,\cdots,t^{\lambda_n}x_n) = t^kg(x_1,\cdots,x_n,t), \text{for}\ \text{some}\ g\in \mathbb C \{x_1,\cdots,x_n,t\}\}.
\end{equation}

Let $\mathfrak a\subset R$ be a nonzero ideal as before, the {\em order} of $\mathfrak a$ with respect to this system of local parameters $T$ and $\lambda\in \mathbb Z^n_{\ge 0}$ is defined as
\begin{equation}\label{ord}
\ord^T_\lambda(\mathfrak a) = \min\{\val^T_{\lambda}(f)|f\in \mathfrak a\}.
\end{equation}

\begin{remark}\label{gen}
If $\mathfrak a$ is generated by $f_1,\cdots,f_r$, then 
\[
\ord^T_\lambda(\mathfrak a) = \min\{\val^T_\lambda(f_i)| i = 1,\cdots,r\}.
\]
Indeed, it is clear that $\min\{\val^T_{\lambda}(f)|f\in \mathfrak a\}\le \min\{\val^T_\lambda(f_i)| i = 1,\cdots,r\}$. On the other hand, if $f\in \mathfrak a$ computes $\ord^T_\lambda(\mathfrak a)$, then we can write $f = \sum_ia_if_i$, for some $a_i\in R$, we have $\val^T_\lambda(f) = \val^T_\lambda(\sum_ia_if_i) \ge \min\{\val^T_\lambda(a_if_i)|i=1,\cdots,r\}\ge \min\{\val^T_\lambda(f_i)|i = 1,\cdots,r\}$.
\end{remark}.

\end{definition}

\begin{definition}\label{Gdef}
Assume $T = \{x_1,\cdots,x_n\}$ is a system of local parameters at $P$ and $\lambda=\{\lambda_1,\cdots,\lambda_n\}\in \mathbb Z^n_{\ge 0}$, we define the {\em slope} function $\mu_P(\lambda,\mathfrak a)$ at $P$ as  
\begin{equation}
\mu_P(\lambda,\mathfrak a) = \frac{\sum_{i = 1}^n \lambda_i}{\ord^T_\lambda(\mathfrak a)}.
\end{equation}
The $T$-{\em stable} {\em rank} of $\mathfrak a$ at $P$ is the infimum of the slope function $\mu_P(\lambda,\mathfrak a)$ with respect to the tuple $\lambda= (\lambda_1,\lambda_2,\cdots,\lambda_n)\in \mathbb Z^n_{\ge 0}$,
\begin{equation}
    \rk^T(P, \mathfrak a) = \inf_\lambda\mu_P(\lambda,\mathfrak a) = \inf_\lambda\frac{\sum_{i = 1}^n \lambda_i}{\ord^T_\lambda(\mathfrak a)}.
\end{equation}
The $G$-{\em stable} {\em rank} of $\mathfrak a$ is defined by taking the infimum of $T$-stable rank with respect to all system of local parameters $T$ at $P$,
\begin{equation}
\rk^G(P, \mathfrak a) = \inf_T(\rk^T(P, \mathfrak a)).
\end{equation}

\end{definition}

If $P\notin V(\mathfrak a)$, we define $\rk^G(P,\mathfrak a) = \infty$, we write $\rk^G(\mathfrak a)$ and $\rk^T(\mathfrak a)$ if $P$ is known in the context. In the following example, we see that an ideal $\mathfrak a$ can have different $T$-stable rank with respect to different system of local parameters $T$ at a point $P$.

\begin{example}
Let $R = \mathbb C[x, y]$, $T = \{x, y\}$ and assume $\mathfrak a = (x^2+2xy+y^2)$ is a principle ideal generated by a polynomial $f(x,y) = x^2+2xy+y^2$, $P=(0,0)$ is the origin. Let $\lambda = (\lambda_1,\lambda_2)\in \mathbb Z^2_{\ge 0}$, then we have
\[
\rk^T(f) = \inf_{\lambda}\frac{\lambda_1+\lambda_2}{\min(2\lambda_1, \lambda_1+\lambda_2,2\lambda_2)} = 1
\]
Let us choose a different system of local parameters $T' = \{u = x + y, v = x - y\}$, then $\mathfrak a = (u^2)$, and $f(u,v) = u^2$, then
\[
\rk^{T'}(f) = \inf_\lambda\frac{\lambda_1+\lambda_2}{2\lambda_1} = \frac{1}{2}
\]
In fact, $\lct_P(u^2) = \frac{1}{2}$, and by Theorem \ref{relation}, we have $\lct_P(\mathfrak a) \le \rk^G(\mathfrak a)$, therefore we get $\rk^G(f) = \frac{1}{2}$.
\end{example}

\begin{example}
Let $\mathfrak a = (x^2y,y^2z,z^2x) \subset \mathbb C[x,y,z]$, $T = \{x, y, x\}$, $P=(0,0, 0)$, then we get
\[
\rk^T(\mathfrak a) = \inf_{\lambda}\frac{\lambda_1 + \lambda_2 + \lambda_3}{\min(2\lambda_1+\lambda_2,2\lambda_2+\lambda_3,2\lambda_3+\lambda_1)} = 1
\]
The ideal $\mathfrak a = (x^2y, y^2z, z^2x)$ is a monomial ideal and we will see later that for a monomial ideal $\mathfrak a$, we have $\rk^G(\mathfrak a) = \lct_P(\mathfrak a)$. Using the fact that $\lct_P(\mathfrak a) = 1$, we obtain $\rk^G(\mathfrak a) = 1$.
\end{example}


%

\begin{remark}\label{ses}
We have a short exact sequence
\begin{equation}
1\to K\to \Aut(\mathbb C\{x_1,\cdots,x_n\})\to \GL(n)\to 1 ,
\end{equation}
where $K$ is a normal subgroup of the group $\Aut(\mathbb C\{x_1,\cdots,x_n\})$ of local holomorphic automorphisms. The morphism $\Aut(\mathbb C\{x_1,\cdots,x_n\})\to \GL(n)$ is given by computing the Jacobian matrix at $(0, \cdots, 0)$. Furthermore, this sequence splits, we have $\Aut(\mathbb C\{x_1\cdots,x_n\}) = K\rtimes \GL(n)$.
\end{remark}

 By Remark \ref{ses}, there is an action of $\GL(n)$ on the set of system of local parameters. Assume  $T=\{x_1,\cdots,x_n\}$ is a system of local parameters at $P$. For $g\in \GL(n)$, $g\cdot T$ is another
 system of local parameters.
We say a system of local parameters $T=\{x_1,\cdots,x_n\}$ is {\em good} for $\mathfrak a$ if $\rk^G(P,\mathfrak a) = \rk^{g\cdot T}(P, \mathfrak a)$
for some $g\in \GL(n)$. In other words, to compute the $G$-stable rank of $\mathfrak a$, it is enough to consider all systems of local parameters obtained from $T$ by actions of $\GL(n)$. 

\begin{example}
Let $f(x,y) = x + y^2 \in \mathbb C[x,y]$, $P = (0, 0)$, we take $T = \{x,y\}$. It can be shown that $\rk^{T}(f) = \frac{3}{2}$. However, if we choose another system of local parameters $T' = \{u = x + y^2, v = y\}$, then $f(u, v) = u$, and $\rk^{T'}(f) = 1$. Indeed, this system of local parameters is optimal, in other words, we can compute the $G$-stable rank in this system of local parameters, and we have $\rk^G(f) = 1$.
\end{example}

\begin{prop}\label{homog}
If $\mathfrak a$ is homogeneous in local parameters $T=\{x_1,\cdots,x_n\}$, then $T$ is good for $\mathfrak a$.
\end{prop}

\begin{proof}
Since $\mathfrak a$ is a homogeneous ideal, we can find a set of generators which are homogeneous polynomials. By Remark \ref{gen}, it is enough to assume that $\mathfrak a$ is generated by a single homogeneous polynomial $f$. Let $g\in K\subseteq \Aut(\mathbb C\{x_1,\cdots,x_n\})$, we can write the action of $g$ on $T$ as following
\[
g(x_i) = x_i + p_i(x_1,\cdots,x_n),
\]
where $p_i\in \mathbb C\{x_1\cdots,x_n\}$ with no constant and degree 1 terms. Since $f$ is a homogeneous polynomial, we have 
\[
\val^T_\lambda(f(g(x_1),\cdots,g(x_n)) \le \val^T_\lambda(f(x_1,\cdots,x_n)).
\]
Let $T'$ be the system of local parameters obtained from $T$ by the action of $g$, then we have 
\begin{equation}\label{comp slopes}
\frac{\sum_{i=1}^n\lambda_i}{\ord^{T'}_\lambda(f)} \ge \frac{\sum_{i=1}^n\lambda_i}{\ord^{T}_\lambda(f)} .
\end{equation}
By Lemma \ref{ses}, given any $h\in \Aut(\mathbb C\{x_1,\cdots,x_n\})$, we can decompose the action of $h$ into an action of $K$ following by an action of $\GL(n)$. By inequality (\ref{comp slopes}), in the system of local parameters obtained by the action of $K$ from $T$, we have larger slope than the slope computed in $T$, hence to compute the $G$-stable rank of $f$, it suffices to consider the action of $\GL(n)$. This shows that $T$ is good for $\mathfrak a$.
\end{proof}

\begin{corollary}
If $f\in \mathbb C[x_1,\cdots,x_n]$ is a homogeneous polynomial of degree $d \ge 2$ and $f$ has isolated singularity at $P = (0,\cdots,0)$. Then $\rk^G(f) = \frac{n}{d}$.
\end{corollary}

\begin{proof}
By Corollary \ref{homog}, the system of local parameters $T = \{x_1,\cdots,x_n\}$ is good for $f$, therefore we only need to consider the group action of $\GL(n)$. 

We claim that $f$ is $\SL(n)$-semistable in the sense of Definition \ref{ss}. Indeed, since $f$ has an isolated singularity at origin, $\frac{\partial f}{\partial x_1},\cdots, \frac{\partial f}{\partial x_n}$ only have a common zero at origin. By \cite{GKZ}, chapter 13, their resultant $\Res(\frac{\partial f}{\partial x_1},\cdots,\frac{\partial f}{\partial x_n})$ is nonzero and invariant under the action of $\SL(n)$. Now assume there is a one parameter subgroup $\lambda: \mathbb G_m\to \SL(n)$, such that 
\[
\lim_{t\to 0} \lambda(t)\cdot v = 0.
\]
Then the resultant $\Res(\frac{\partial f}{\partial x_1},\cdots,\frac{\partial f}{\partial x_n})$ is 0, which is impossible. This proves the claim.

By the claim that $f$ is $\SL(n)$-semistable, the corollary follows immediately from Proposition 2.6 in \cite{G-stable rank}
\end{proof}

\subsection{Relation to log canonical threshold}
Let $X = \text{Spec}(R)$ be a nonsingular irreducible complex affine variety, $\mathfrak a\subset R$ is an ideal, $P\in V(\mathfrak a)$. If $\mathfrak a = (f_1,\cdots,f_r)\subseteq R$ is a nonzero ideal, recall the Definition \ref{lct}, the log canonical threshold $\lct_P(\mathfrak a)$ of the ideal $\mathfrak a$ at point $P$ is

\begin{equation}
    \lct_P(\mathfrak a) = \sup\{s>0|\frac{1}{(\sum_{i = 0}^r|f_i|^2)^s}  \
    \text{is}\ \text{integrable}\ \text{around}\ P\}.
\end{equation}

Theorem \ref{relation} says that the log canonical threshold of an ideal is less than or equal to the $G$-stable rank of that ideal
\begin{equation}
\lct_P(\mathfrak a) \le \rk^G(P,\mathfrak a).
\end{equation}

\begin{proof}[Proof of Theorem \ref{relation}]
Let $s>0$ be such that $\frac{1}{(\sum_{i = 0}^r|f_i|^2)^s}$ is integrable around $P$, then there is a neighborhood $U_P$ of $P$, such that

\[
\int_{U_P} \frac{dV}{(\sum_{i = 0}^r|f_i|^2)^s} < C < \infty,
\]
for some constant $C$. Choose a system of local parameters $T = \{x_1,\cdots,x_n\}$ at $P$, let $\lambda = (\lambda_1,\cdots,\lambda_n)\in\mathbb Z^n_{\ge0}$, then let $t\in \mathbb C^*$ act on the coordinates by $x_i\to t^{\lambda_i}x_i$. We denote $t\cdot U_P$ for the image of $U_P$ under this action. If $|t|<1$, we have $t\cdot U_P\subset U_P$, therefore
\begin{equation}
\int_{t\cdot U_P} \frac{dV}{(\sum_{i = 0}^r|f_i|^2)^s} < \int_{U_P} \frac{dV}{(\sum_{i = 0}^r|f_i|^2)^s} <C
\end{equation}

Let $y_i = t^{-\lambda_i} x_i$ for $i = 1,\cdots, n$, then we have
\begin{equation}\label{int}
    \int_{t\cdot U_P} \frac{dx_1d\bar{x}_1\cdots dx_nd\bar x_n}{(\sum_{i = 0}^r|f_i(x_1,\cdots,x_n)|^2)^s}
    = \int_{U_P}\frac{|t|^{2\sum_{i=1}^n \lambda_i}dy_1d\bar y_1\cdots dy_nd\bar y_n}{(\sum_{i=1}^r|f_i(t^{\lambda_1}y_1,\cdots,t^{\lambda_n}y_n)|^2)^s} < C
\end{equation}

Recall the Definition \ref{val&ord} for the valuation (\ref{val}) and order of an ideal (\ref{ord}), we can write
\[
\sum_{i=0}^r|f_i(t^{\lambda_1}y_1,\cdots,t^{\lambda_n}y_n)|^2= |t|^{2\min_i(\val^T_\lambda(f_i))}(\sum_{i=1}^r|\tilde f_i(y_1,\cdots,y_n, t)|^2) 
\]
\[= |t|^{2\ord_\lambda^T(\mathfrak a)}(\sum_{i=1}^r|\tilde f_i(y_1,\cdots, y_n, t)|^2),
\]
for some $\tilde f_i(y_1,\cdots, y_n, t) \in \mathbb C\{y_1,y_2,\cdots,y_n,t\}$. In particular, we know that $\sum_{i=1}^r|\tilde f_i(y_1,\cdots, y_n, 0)|^2$ is not constantly zero in $U_P$. So we can find a point $Q\in U_P$ such that $ 0<\sum_{i=1}^r|\tilde f_i(Q, 0)|^2 < B$ for some constant $B > 0$. By the continuity, there is a neighborhood $U_Q$ such that $Q\in U_Q \subset U_P$ and some $\epsilon > 0$, such that $ 0<\sum_{i=1}^r|\tilde f_i(y_1,\cdots,y_n, t)|^2 < B$ for any $(y_1,\cdots,y_n) \in U_Q$ and $0 \le t < \epsilon$. We can write integral (\ref{int}) as

\begin{equation}\label{new_int}
    \int_{U_P}\frac{|t|^{2\sum_{i=1}^n \lambda_i}dy_1d\bar y_1\cdots dy_nd\bar y_n}{|t|^{2\ord_\lambda^T(\mathfrak a)}(\sum_{i=1}^r|\tilde f_i(y_1,\cdots, y_n, t)|^2)^s} = \int_{U_P}\frac{|t|^{2(\sum_{i=1}^n \lambda_i - s\ord_\lambda^T(\mathfrak a))}dV}{(\sum_{i=1}^r|\tilde f_i(y_1,\cdots, y_n, t)|^2)^s}< C.
\end{equation}

Therefore we get
\[
|t|^{2(\sum_{i=1}^n \lambda_i - s\ord_\lambda^T(\mathfrak a))}\int_{U_Q}\frac{dV}{B^s}
=\int_{U_Q}\frac{|t|^{2(\sum_{i=1}^n \lambda_i - s\ord_\lambda^T(\mathfrak a))}dV}{B^s}
<\int_{U_P}\frac{|t|^{2(\sum_{i=1}^n \lambda_i - s\ord_\lambda^T(\mathfrak a))}dV}{(\sum_{i=1}^r|\tilde f_i(y_1,\cdots, y_n, t)|^2)^s} < C
\]
for all $t\in [0,\epsilon)$, where $\epsilon > 0$. Since $\int_{U_Q}\frac{dV}{B^s} > 0$, we must have 

\[
\sum_{i=1}^n\lambda_i-s\ord^T_{\lambda}(\mathfrak a) \ge 0.
\]
Hence we get $s\le \frac{\sum_{i=1}^n\lambda_i}{\ord^T_{\lambda}(\mathfrak a)}$. This holds for any system of local coordinates and $\lambda = (\lambda_1,\cdots, \lambda_n)\in\mathbb Z^n_{\ge0}$, therefore
\[
\lct_P(\mathfrak a) \le \rk^G(P, \mathfrak a).
\]

This completes the proof of Theorem \ref{relation}.
\end{proof}
\begin{example}
Let $f=x_1^{u_1}+x_2^{u_2}+\cdots +x_n^{u_n}\in \mathbb  C[x_1,\cdots,x_n]$, $P = (0,\cdots,0)$ be the origin. It was shown in \cite{Mus} that $\lct_P(f) = \min(1,\sum_{i=1}^n\frac{1}{u_i})$. However, we will show that $\rk^G(f) =\sum_{i=1}^n \frac{1}{u_i}$. So we get $\lct_P(f)\le \rk^G(f)$. If $n = 3$, $u_1 = u_2 = u_3 = 2$, then we have $f = x_1^2+x_2^2+x_3^2$ and $\lct_P(f) = 1 < \rk^G(f) = \frac{3}{2}$.
\end{example}

Suppose $\mathfrak a = (m_1,\cdots,m_r) \subset \mathbb C[x_1,\cdots,x_n]$ is a proper nonzero ideal generated by monomials $\{m_1,\cdots, m_r\}$ and let $P = (0,\cdots, 0)$ be the origin. Given $u = (u_1,\cdots, u_n)\in \mathbb Z^n_{\ge 0}$, we write $x^u = x_1^{u_1}\cdots x_n^{u_n}$. The $Newton$ $Polyhedron$ of $\mathfrak a$ is 
\begin{equation}
P(\mathfrak a) =\text{convex}\ \text{hull}\ (\{u\in\mathbb Z^n_{\ge 0} |x^u\in \mathfrak a\}).
\end{equation}

It was shown in \cite{Mus} that
\begin{equation}
\lct_P(\mathfrak a) = \max\{\nu\in\mathbb R_{\ge 0}|(1,1,\cdots,1)\in \nu\cdot P(\mathfrak a)\}.
\end{equation}
In other words, $\lct_P(\mathfrak a)$ is equal to the largest $\nu$ such that $\sum_{i=1}^n\lambda_i \ge \nu\cdot \min_{u\in P(\mathfrak a)}\langle u,\lambda\rangle$ for any $\lambda = (\lambda_1,\cdots,\lambda_n)\in \mathbb Z^n_{\ge 0}$, where we use the standard inner product $\langle u,\lambda\rangle = \sum_{i=1}^n u_i\lambda_i$.

Theorem \ref{equal} says that the log canonical threshold is equal to the $G$-stable rank for monomial ideals. More precisely, suppose $\mathfrak a\subset \mathbb C[x_1,\cdots,x_n]$ is a proper nonzero ideal generated by monomials and $P=(0,\cdots,0)$ is the origin. Then
\begin{equation}
\lct_P(\mathfrak a) = \rk^G(\mathfrak a).
\end{equation}

\begin{proof}[Proof of Theorem \ref{equal}]
Let $\mathfrak a = (m_1,\cdots,m_r)$ and $\{m_i = x_1^{l_{i1}}x_2^{l_{i2}}\cdots x_n^{l_{in}}\}_{i=1,\cdots,r}$ be a set of generators, where $l_i = (l_{i1},\cdots,l_{in})\in \mathbb Z^n_{\ge 0}$, we have 
\[
\lct_P(\mathfrak a) = \max\{\nu\in\mathbb R_{\ge 0}|(1,\cdots,1)\in \nu \cdot P(\mathfrak a)\}
\]
\[
=\max\{\nu\in \mathbb R_{\ge 0}| \sum_j\lambda_j \ge \nu\cdot \min_{u\in P(\mathfrak a)}
\langle u,\lambda\rangle, \forall \lambda\in \mathbb Z^n_{\ge 0}\}
\]
\[
=\max\{\nu\in \mathbb R_{\ge 0}| \sum_j\lambda_j \ge \nu\cdot \min_{i}
\langle l_i,\lambda\rangle, \forall \lambda\in \mathbb Z^n_{\ge 0}\}
\]
\[
=\max\{\nu\in \mathbb R_{\ge 0}| \nu \le \frac{\sum_j\lambda_j}{\min_i(\sum_j l_{ij}\lambda_j)}, \forall \lambda\in \mathbb Z^n_{\ge 0}\}.
\]
In this system of local parameters $T = \{x_1,\cdots,x_n\}$, we have $\ord^T_\lambda(\mathfrak a) = \min_i(\sum_j l_{ij}\lambda_j)$, therefore we get
\[
\lct_P(\mathfrak a) = \max\{\nu\in \mathbb R_{\ge 0}| \nu \le \frac{\sum_i\lambda_i}{\ord^T_\lambda(\mathfrak a)}, \forall \lambda\in \mathbb Z^n_{\ge 0}\}
\]
\[
= \max\{\nu\in \mathbb R_{\ge 0}| \nu \le \rk^T(\mathfrak a)\} = \rk^T(\mathfrak a).
\]
Since log canonical threshold does not depend on the local coordinates, hence we have 
\[
\lct_P(\mathfrak a) = \rk^G(\mathfrak a).
\]
This completes the proof of Theorem \ref{equal}.
\end{proof}
\begin{example}
Suppose $\mathfrak a = (x_1^{u_1},\cdots,x_n^{u_n})$ and $P = (0,\cdots,0) $, then $\lct_P(\mathfrak a) = \sum_{i=1}^n\frac{1}{u_i}$. Therefore we also have $\rk^G(\mathfrak a) = \sum_{i=1}^n\frac{1}{u_i}$.
\end{example}

\subsection{Some properties of G-stable rank for ideals}

Some results  for log canonical threshold can be found in~\cite{Mus}. Here, we also prove similar results for the $G$-stable rank of ideals. 

\begin{prop}
If $\mathfrak a\subseteq \mathfrak b$ are nonzero ideals on $X$, then we have $\lct_P(\mathfrak a)\le \lct_P(\mathfrak b)$ and $\rk^G(P,\mathfrak a)\le \rk^G(P,\mathfrak b)$. 
\end{prop}
\begin{proof}
The first inequality was shown in \cite{Mus}. If $P\notin V(\mathfrak b)$, it is trivial. Assume $P\in V(\mathfrak b)$, let $T$ be a system of local parameters at $P$ and $\lambda=(\lambda_1,\cdots,\lambda_n)\in \mathbb Z^n_{\ge 0}$. Since $\mathfrak a\subseteq \mathfrak b$, we have $\ord^T_\lambda(\mathfrak b) \le  \ord^T_\lambda(\mathfrak a)$, therefore $\mu_P(\lambda,\mathfrak a) \le \mu_P(\lambda,\mathfrak b)$, it follows immediately that $\rk^G(P,\mathfrak a)\le \rk^G(P,\mathfrak b)$.
\end{proof}

\begin{prop}
We have $\lct_P(\mathfrak a^r) = \frac{\lct_P(\mathfrak a)}{r}$ and $\rk^G(P,\mathfrak a^r) = \frac{\rk^G(P,\mathfrak a)}{r}$ for every $r\ge 1$.
\end{prop}
\begin{proof}
The first claim was shown in \cite{Mus} and the second claim follows from the fact that $\ord^T_\lambda(\mathfrak a^r) = r\cdot\ord^T_\lambda(\mathfrak a)$ and Definition \ref{Gdef}.
\end{proof}

\begin{prop}
If $\mathfrak a$ and $\mathfrak b$ are ideals of $X$, then 
\[
\frac{1}{\lct_P(\mathfrak a\cdot\mathfrak b)} \le \frac{1}{\lct_P(\mathfrak a)} + \frac{1}{\lct_P(\mathfrak b)}, \quad
\frac{1}{\rk^G(P,\mathfrak a\cdot\mathfrak b)}\le \frac{1}{\rk^G(P,\mathfrak a)}+\frac{1}{\rk^G(P,\mathfrak b)}.
\]
\end{prop}
\begin{proof}
The first inequality was shown in \cite{Mus}, we show the second inequality. Let $T$ be a system of local parameters at $P$ and $\lambda = (\lambda_1,\cdots,\lambda_n)\in \mathbb Z^n_{\ge 0}$, then  
\[
\ord_\lambda^T(\mathfrak a\cdot\mathfrak b) = \ord_\lambda^T(\mathfrak a) +\ord_\lambda^T(\mathfrak b).
\]
Therefore, we have
\[
\sup_{T,\lambda}\frac{\ord_\lambda^T(\mathfrak a\cdot\mathfrak b)}{\sum\lambda_i}= \sup_{T,\lambda}\left(\frac{\ord_\lambda^T(\mathfrak a)}{\sum\lambda_i}+\frac{\ord_\lambda^T(\mathfrak b)}{\sum\lambda_i}\right)
\]
\[
\le \sup_{T,\lambda}\frac{\ord_\lambda^T(\mathfrak a)}{\sum\lambda_i}+\sup_{T,\lambda}\frac{\ord_\lambda^T(\mathfrak b)}{\sum\lambda_i}.
\]
\end{proof}

The following two propositions are from \cite{Mus}. A lot of evidence suggests the same results for $G$-stable rank of ideals, we give them as conjectures. 
\begin{prop}
If $H\subset X$ is a nonsingular hypersurface such that $\mathfrak a\cdot \mathcal O_H$ is nonzero, then $\lct_P(\mathfrak a\cdot \mathcal O_H)\le \lct_P(\mathfrak a)$.
\end{prop}

\begin{prop}
If $\mathfrak a$ and $\mathfrak b$ are ideals on $X$, then
\[
\lct_P(\mathfrak a + \mathfrak b) \le \lct_P(\mathfrak a) + \lct_P(\mathfrak b)
\]
for every $P\in X$.
\end{prop}

\begin{conjecture}
If $H\subset X$ is a nonsingular hypersurface such that $\mathfrak a\cdot \mathcal O_H$ is nonzero, then $\rk^G(P,\mathfrak a\cdot\mathcal O_H)\le \rk^G(P,\mathfrak a)$ for any $P\in H$.
\end{conjecture}

\begin{conjecture}
Let $\mathfrak a$ and $\mathfrak b$ be two nonzero proper ideals of $X$, then for any point $P\in X$, we have 
\[
    \rk^G(P, \mathfrak a+\mathfrak b) \le \ \rk^G(P, \mathfrak a) +  \rk^G(P, \mathfrak b).
\]
\end{conjecture}

\section{Acknowledgements}
I would like to thank Harm Derksen for direction and discussion, and Visu Makam for comments on an earlier draft of this paper.

\end{document}